\documentclass[a4paper,12pt]{article}

\newtheorem{algorithm}{Algorithm}

\newtheorem{remark}{Remark}

\usepackage{booktabs}

\usepackage{algorithm}
\usepackage[noend]{algpseudocode}
\makeatletter
\def\BState{\State\hskip-\ALG@thistlm}
\makeatother
\usepackage{algorithmicx}

\usepackage{amsmath}
\usepackage{amssymb}
\usepackage{amsfonts}
\usepackage[dvipdfmx]{graphicx}
\usepackage{ dsfont }

\usepackage[normalem]{ulem}
\usepackage[dvipsnames]{xcolor}

\usepackage{tikz}
\usepackage{lipsum}

\newcommand{\bs}{\boldsymbol}
\newcommand{\mr}[1]{\mathrm{#1}}

\newcommand{\I}{\mathcal{I}}

\newcommand{\Rb}{\mathbb{R}}

\title{\bf Power grid transient stabilization using Koopman model predictive control}

\begin{document}

\author{Milan Korda$^1$, Yoshihiko Susuki$^2$, Igor Mezi{\'c}$^1$}

\footnotetext[1]{Milan Korda and Igor Mezi{\'c} are with the University of California, Santa Barbara,\; {\tt milan.korda@engineering.ucsb.edu, mezic@engineering.ucsb.edu}
Yoshihiko Susuki is with Osaka Prefecture University,\; {\tt susuki@eis.osakafu-u.ac.jp}}

\date{Draft of \today}

\maketitle

\begin{abstract}
This work addresses the problem of transient stabilization of a power grid, following a destabilizing disturbance. The model considered is the cascade interconnection of seven New England test models with the disturbance (e.g., a powerline failure) occurring in the first grid and propagating forward, emulating a wide-area blackout. We consider a data-driven control framework based on the Koopman operator theory, where a linear predictor, evolving on a higher dimensional (embedded) state-space, is built from observed data and subsequently used within a model predictive control (MPC) framework, allowing for the use of efficient computational tools of linear MPC to control this highly nonlinear dynamical system.
\end{abstract}

\begin{flushleft}\small
\end{flushleft}

\section{Introduction}

Transient stabilization is of vital importance for emergency control of large-scale interconnected power grids.  
The stabilization problem is associated with the control of electromechanical dynamics of coupled synchronous generators when subjected to a large disturbance;  see e.g., \cite{kundur1994}.  
Since it is typical of nonlinear and large-scale control problems, it is a well-established subject with a long history of research: see e.g. \cite{bazanella1999,galaz2003,ortega2005,cornelius2013}.  
Failure of transient stabilization is recognized as one cause of large blackouts such as the September 2003 blackout in Italy; see \cite{corsi2004} for details.   
All of the transient stabilization controllers proposed in the literature are model-driven, to the best of the authors' knowledge. 
In addition  to 
the inherent nonlinear and large-scale nature, 
the wide-spread introduction of renewable energy resources with uncertain characteristics 
makes it extremely challenging to derive an accurate deterministic model of the power grid. 
Therefore, an alternative, data-driven 
method for the design of transient stabilization controllers is currently required.


The stabilization controller proposed in this paper is based on the Koopman operator model predictive control (MPC) of~\cite{korda2016linear}, where a \emph{linear} predictor is constructed from observed data generated by the \emph{nonlinear} dynamical system. The distinguishing feature of the predictor is the fact that its state evolves on a higher-dimensional, embedded, state-space, thereby being able to capture the nonlinear behavior of the underlying dynamics. This predictor is subsequently used within a linear model predictive control 
scheme, thereby allowing for the use of highly efficient linear MPC tools to control the nonlinear dynamical system. Importantly, the dimension of the embedding space does not affect the computational complexity of the optimization problem solved by the MPC and therefore the control input can be evaluated very fast, allowing for a real-time deployment.

\begin{figure}[t]
	\begin{picture}(50,200)
	\put(80,0){\includegraphics[width=120mm]{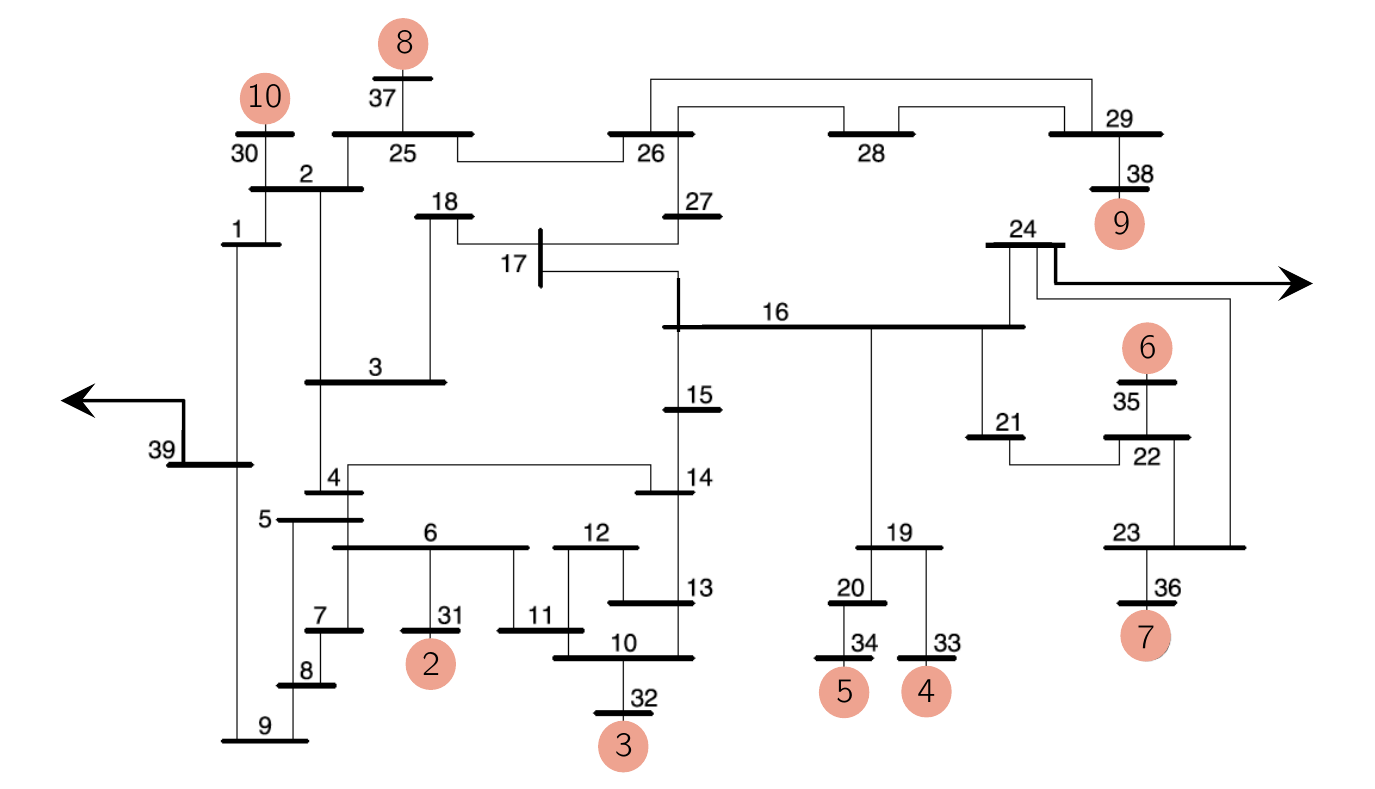}}
	\end{picture}
	\caption{\footnotesize One-line diagram of single New England test model within the cascade. The first generator, normally connected to bus 39, is replaced by a connection to the previous grid of the cascade. Generator one in the first grid represents the infinite bus.}
	\label{fig:NEoneGrid}
\end{figure}

In this paper, we numerically demonstrate the Koopman MPC for transient stabilization of power grids. 
Numerical simulations are conducted for the cascade interconnection of seven New England test models shown in \cite{susuki2012coherent} (see Figures~\ref{fig:NEcascade} and \ref{fig:NEoneGrid}). The cascade connection was invented for exploring the mechanism of disturbance propagation based on the theory of coherent swing instability in \cite{susuki2011coherent}.    
Transient stabilization of the cascaded grid is investigated using the data-driven methods of~\cite{korda2016linear}; this is the first demonstration of these methods for control of power grid dynamics. The results are promising, achieving a successful stabilization of the cascaded grid without the model knowledge, with a distributed control structure (one controller per grid) and fast computation time; future work will investigate and compare the efficacy of the proposed method on different power grid models and control setups.


\section{Problem statement}
\begin{figure}[t]
	\begin{picture}(50,40)
	\put(0,0){\includegraphics[width=160mm]{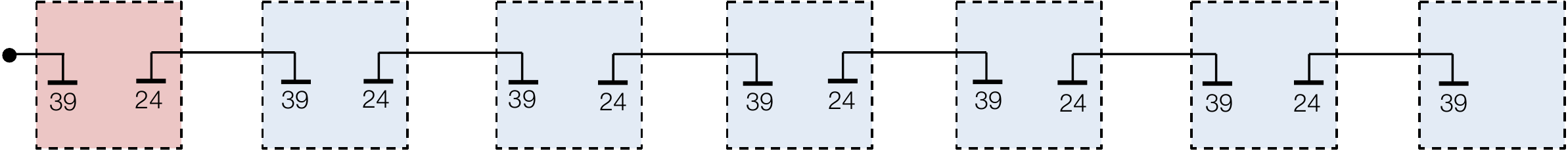}}
	\end{picture}
	\caption{\footnotesize Cascade interconnection of seven New England test models. The fault occurs in the first grid near bus 39 (modeled by adding a small impedance between bus 39 and the ground).}
	\label{fig:NEcascade}
\end{figure}

We use the so-called nonlinear swing equations (see e.g. \cite{kundur1994}) for modeling and analysis of coupled swing dynamics in the cascaded system of seven New England test models. The short-term electromechanical dynamics of generator $j$ in unit grid \#$i$ $(j=2,\ldots,10, i=1,\ldots,7)$ in Figure~\ref{fig:NEcascade} are represented by the so-called swing equations as follows:
\begin{equation}
\left.
\begin{array}{rcl}
\displaystyle \frac{d\delta_{ij}}{d t}\, &=& \omega_{ij},
\\\noalign{\vskip +2.0mm}
\displaystyle \frac{H_j}{\pi f_{\rm b}}\frac{d\omega_{ij}}{d t} &=& P_{{\rm m}\,ij}(1+u_{ij}) - D_j\omega_{ij}
\\\noalign{\vskip +1.0mm}
& &\hspace{-20mm} -V_{ij}V_{11}\{G_{ij,11}\cos(\delta_{ij}-\delta_{11})+B_{ij,11}\sin(\delta_{ij}-\delta_{11})\}
\\\noalign{\vskip +2.0mm}
& & \multicolumn{1}{r}{\hspace{-20mm}\displaystyle -V^2_{ij}G_{ij}-\sum^{7}_{k=1,k\neq i}\sum^{10}_{l=2,l\neq j}V_{ij}V_{kl}\{G_{ij,kl}\cos(\delta_{ij}-\delta_{kl})}
\\\noalign{\vskip +2.0mm}
& & 
\multicolumn{1}{r}{\hspace{-10mm} +B_{ij,kl}\sin(\delta_{ij}-\delta_{kl})\}.}
\end{array}
\makebox[+0.5em]{}\right\}
\label{eqn-SSK:SE_7NE}
\end{equation}
The variable $\delta_{ij}$ represents the angular position of rotor in generator $j$ in unit grid \#$i$ with respect to the infinite bus (see e.g. \cite{kundur1994}) and is in radians [rad].  
The variable $\omega_{ij}$ represents the deviation of rotor speed in generator $i$ relative to system angular frequency $2\pi f_{\rm b}$ and is in radians per second [rad/s].  
The variable $\delta_{11}$ is the angular position of the infinite bus and becomes constant from its definition.  
The parameters $H_j$, $P_{{\rm m}\,ij}$, $D_j$, $V_{ij}$, $G_{ij,ij}$, and $G_{ij,kl}+{\rm i}B_{ij,kl}$ are in per unit system except for $H_j$ and $D_j$ in seconds [sec].  
The parameter $H_j$ denotes the per-unit time inertia constant of generator $j$, $D_j$ denotes its dampling coefficient.  
The constant $P_{{\rm m}\,ij}$ is the nominal mechanical input power to generator $j$ in unit grid \#$i$, and $V_{ij}$ is the internal voltage of generator $j$ in unit grid \#$i$.  
They are assumed to be constant.  
The constant $G_{ij,ij}$ denotes the internal conductance of generator $j$ in unit grid \#$i$, and $G_{ij,kl}+{\rm i}B_{ij,kl}$ denotes the transfer admittance between generators $j$ in unit grid \#$i$ and $\l$ in unit grid \#$k$.  
The constant $V_{11}$ is the voltage of the infinite bus, and $G_{ij,11}+{\rm i}B_{ij,11}$ is the transfer admittance between generator $j$ in unit grid \#$i$ and the infinite bus.  
The impedance $G_{ij,kl}+{\rm i}B_{ij,kl}$ is the parameters that change as the network topology changes. 
The control input $u_{ij}$ is the mechanical input power expressed as the fraction of the nominal value $P_{\mr{m}\,ij}$.

The fault causing the instability is modeled by adding a small impedance, $10^{-7}\,\Omega$, between bus 39 and the ground. The fault occurs in grid $\#1$ at time $t_{\rm f}  = 0.87\,\mr{s}$ and the line $1-39$ trips (i.e., the this line is removed from the model) at time $t = 1\,\mr{s}$.  Prior to the fault time $t_{\rm f}$, the system operates in a steady-state condition given by a power-flow computation.

The goal is to control the instability by adjusting the generator mechanical input power\footnote{In this work, we chose to control the mechanical input power. Alternatively, one could also control the generator voltages $V_{ij}$.} while observing the following engineering requirements:
\begin{itemize}
\item Frequency deviations convergence to zero.  
This implies for the current power grid model that all of the generators settle down to the nominal rotating frequency (back to a state of frequency synchronization).
\item Maximum frequency deviation does not exceed a given bound.  
This is a mandatory requirement in the practical power grid; otherwise, in order to avoid a serous damage of turbine blades in a power plant, the corresponding power plant (generator) will be removed from the grid by a protective relay.
\item Control action is within given bounds.
This requirement is also practical because such controllers contain a limiter (or saturation device) in order to prevent exceeding physical limitations of the generators.
\end{itemize}

\section{Control strategy}
\begin{figure*}[t]
	\begin{picture}(60,120)
	\put(65,0){\includegraphics[width=140mm]{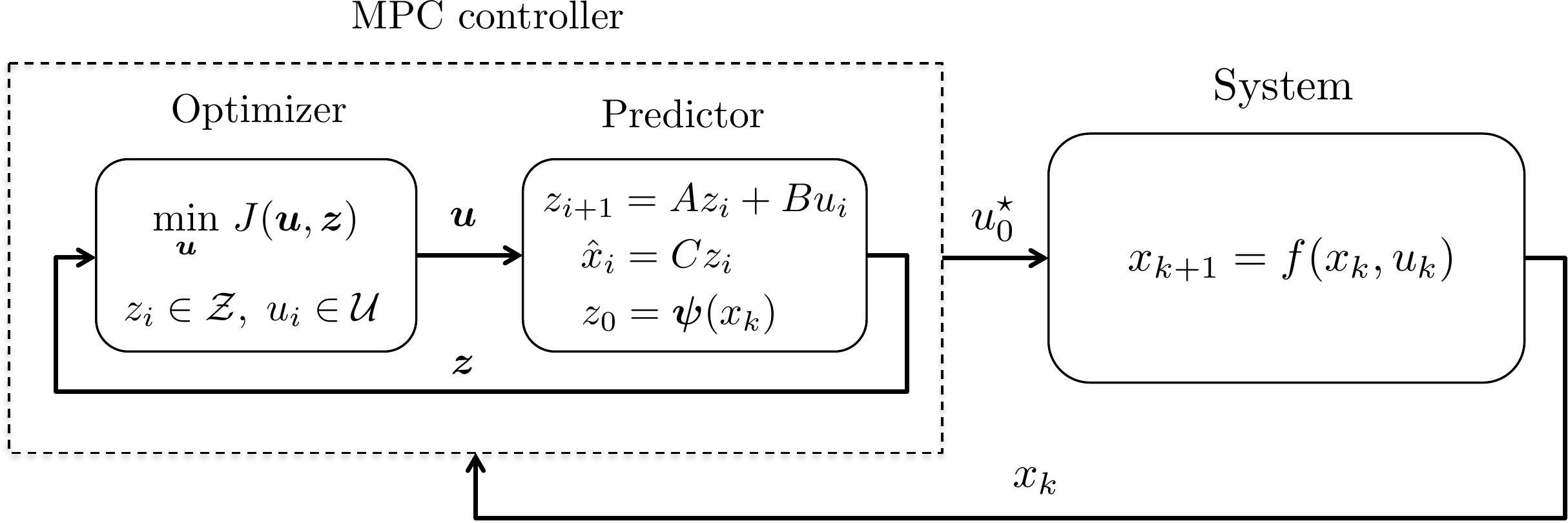}}
	\end{picture}
	\caption{\footnotesize Conceptual depiction of the Koopman MPC, where $x \in \Rb^n$ is the true state of the dynamical system, $z \in \Rb^N$ is the embedded state, $\bs \psi:\mathbb{R}^n \to \Rb^N$ is a nonlinear embedding mapping, $\hat x$ is prediction of the true state using the linear predictor and $u$ is the control input. Bold symbols denote vectorization along the prediction horizon $N_p$, i.e., $\bs u = (u_0,\ldots,u_{N_p-1})$, $\bs z = (z_0,\ldots,z_{N_p})$. The symbol $u_0^\star$ denotes the first component of the optimal solution $\bs u^\star = (u_0^\star,\ldots, u_{N_p-1}^\star)$ to the optimization problem solved by the optimizer. This optimization problem is solved repeatedly at each time step $k$ of the closed-loop operation, with the state of the predictor re-initialized to $z_0 = \bs\psi(x_k)$. Note the difference between the ``physical'' time~$k$ and the time of the predictor $i$. The objective function $J$ is a convex quadratic function and the constraint sets $\mathcal{Z}$ and $\mathcal{U}$ are polyhedra. The optimization problem is therefore a convex quadratic program.}
	\label{fig:MPC}
\end{figure*}

In order to control the instability we use the Koopman model predictive control proposed in~\cite{korda2016linear}. The conceptual scheme of the strategy is depicted and described in Figure~\ref{fig:MPC}. The main components of the MPC controller are a \emph{predictor} and an \emph{optimizer}. The distinguishing feature of the Koopman MPC controller is the fact that the predictor is in the form of a \emph{linear} dynamical system evolving on an embedded (or lifted) state space of larger dimension than the dimension of original state space. Contrary to classical local linearizion techniques, this predictor is valid \emph{globally} (or in a large subset of the state space) and fixed once and for all. The linearity of the predictor and the freedom in the choice of the embedding mapping implies that the optimization problem solved can be rendered \emph{convex} quadratic, even if the original problem had non-convex objective function and constraints and nonlinear dynamics (see~\cite{korda2016linear} for details and \cite{korda2017convergence} for a theoretical analysis of such predictors).  Importantly, by tranforming the optimization problem in the so-called dense form (see Section~\ref{sec:solvingMPC}), the computational complexity can be rendered independent of the dimension of the embedded state-space; see Remark~\ref{rem:denseForm} and \cite{korda2016linear} for details.

\begin{remark} We note that the depiction of the MPC controller in Figure~\ref{fig:MPC} is only conceptual and in real applications, the dynamics of the predictor is a part of the constraints of the optimization problem solved by the optimizer.
\end{remark}

\subsection{Koopman MPC for power grid}\label{sec:KoopPG}
Now we describe the use of the general scheme from Figure~\ref{fig:MPC} for the cascade instability control. In this case, the state is given by $x = (\delta,\omega)$. We choose the embedding mapping $\bs \psi $ to be
\begin{equation}\label{eq:embedMap}
 \bs\psi(x) = \begin{bmatrix}\cos\delta \\ \sin\delta \\ \omega\end{bmatrix}
\end{equation} 
and for the time being we assume that a predictor of the form
\begin{align}\label{eq:pred_lin}
z_{i+1} & = Az_i + B u_i, \\
\hat{x}_i &= C z_i \nonumber, 
\end{align}
is given; here, $u \in \mathbb{R}^m$ is the control input with $m$ being the number of manipulated variables\footnote{In this work we shall consider distributed Koopman MPC with one controller designed per grid (see Section~\ref{sec:distMPC}). In this case, $m = 9$ is the number of generators in each grid. Of course, other control configurations can be considered, leading to different dimensions of the manipulated variable $u$. The discussion of Section~\ref{sec:KoopPG} is fully general, applying to any control configuration.}; see Section~\ref{sec:predConst} for a method to construct such predictor from data. The Koopman MPC controller then solves at each time step $k$ of the closed-loop operation the optimization problem
\begin{equation}\label{opt:MPC}
\begin{array}{ll}
\underset{\textstyle u_0,\ldots, \textstyle u_{N_p-1} }{\mbox{minimize}} & \sum_{i=0}^{N_p-1} z_i^\top Q z_i + u_i^\top R u_i  \\
\mbox{subject to} & z_{i+1} = Az_i + Bu_i, \quad  \hspace{5mm} i = 0,\ldots, N_p-1 \\
& u_{\mr{min}} \le u_i \le u_{\mr{max}},  
\hspace{1mm}\quad \quad i = 0,\ldots, N_p-1\\
& z_{\mr{min}} \le z_i \le z_{\mr{max}}, \hspace{2mm}\quad \quad i = 0,\ldots, N_p-1\\
\mbox{parameter} & z_0 = {\bs \psi}(x_k) ,
\end{array}
\end{equation}
where $\cdot^\top$ denotes the transpose of a vector and $Q\succeq 0$, $R\succeq 0$ are given symmetric positive semidefinite weighting matrices. These matrices are chosen according to the desired control goals; for example, a natural goal is to minimize the sum of square deviations of the angular frequencies from the synchronous frequency $2\pi 60 \,\mr{rad}/\mr{s}$, i.e., to minimize $\sum_{i=0}^{N_p-1}\|\omega_i\|_2^2$, where $\|.\|_2$ denotes the Euclidean norm. Given the choice of the embedded state $z = \bs\psi(x)$ in~(\ref{eq:embedMap}), this can be expressed by choosing
\[
Q = \mr{bdiag}(0_{ n_{\mr{gen}} \times n_{\mr{gen}} }, 0_{n_{\mr{gen}}\times n_{\mr{gen}}}, I_{n_{\mr{gen}}}  ), \quad R = 0,
\]
where $n_{\mr{gen}}$ is the number of generators (in each grid in the case of distributed control or the total number of generators in the case of centralized control) and $0_{n\times n}$ and $I_n$ denote, respectively, the zero and identity matrices of size $n\times n$ and $\mr{bdiag}(\cdot,\ldots,\cdot)$ denotes a block diagonal matrix composed of the arguments. Similarly, a constraint that the maximum deviation of the generator angles is no more than $\theta_{\mr{max}}$ and maximum deviation of the  angular frequencies no more than $\omega_{\mr{max}}$ is imposed by choosing
\[
z_\mr{max} = \begin{bmatrix} \cos(\theta_{\mr{max}}) \\ \sin( \theta_{\mr{max}}) \\ \omega_{\mr{max}}  \end{bmatrix},\quad z_{\mr{min}} = - z_{\mr{max}}.
\]

The Koopman MPC algorithm for the power grid is summarized in Algorithm~\ref{alg:MPC}.

\begin{remark}\label{rem:denseForm}
Note that, strictly speaking, the minimization in~(\ref{opt:MPC}) should be over both $u_i$'s and $z_i$'s. However, since $z_0$ and $u_0,\ldots,u_{N-1}$ uniquely determine $z_0,\ldots, z_N$ via $z_{i+1} = Az_i + B u_i$, the $z_i$'s can be eliminated (solved for), obtaining the so-called dense form of MPC. This eliminates the dependence on $z_i$'s and therefore on the dimension of the embedded state $z$; see~\cite{korda2016linear} for details of this transformation.
\end{remark}

\begin{algorithm}[t]
\caption{Koopman MPC -- closed-loop operation}\label{alg:MPC}
\begin{algorithmic}[1]
\Require  Predictor ($A$,$B$), Cost matrices $(Q,R)$, bounds $(u_{\mr{max}},\,u_{\mr{min}},\,z_{\mr{max}},\,z_{\mr{min}})$.
\For{$k=0,1,\ldots$}
\State Measure $x_k = (\delta_k,\omega_k) $ on the real system
\State Set $z_0 = \bs\psi(x_k)=  [\cos\delta_k^\top, \sin\delta_k^\top,\omega_k^\top]^\top$
\State		Solve~(\ref{opt:MPC}) to get an optimal solution $(u_0^\star,\ldots,u_{N_p-1}^\star)$
\State Apply $u_k := u_0^\star$ to the real system
\EndFor
\end{algorithmic}
\end{algorithm}

\subsection{Solving optimization problem (\ref{opt:MPC})}\label{sec:solvingMPC}
The optimization problem~(\ref{opt:MPC}) is a convex quadratic program and therefore can be efficiently solved by a number of off-the-shelf solvers (e.g., CPLEX, GUROBI, MOSEK etc). In addition,  the special structure of the MPC problem can be exploited by tailored algorithms. In our case, we first eliminate $z_i$'s in~(\ref{opt:MPC}), arriving at the so-called dense form of the problem
\begin{equation}\label{eq:MPC_dense}
\begin{array}{ll}
\underset{U \in \Rb^{mN_p}}{\mbox{minimize}} & U^\top H U^\top  + z_0^\top G U  \\
\mbox{subject to} &  L U + M z_0 \le c \\
\mbox{parameter} & z_0 = {\bs \psi}(x_k) ,
\end{array}
\end{equation}
where the decision variable $U$ is
\[
U = [u_0^\top,\ldots, u_{N_p-1}^\top]^\top
\]
and the data matrices are given by
 \[
H = {\mathbf{R}} +{ \mathbf{B}}^\top\mathbf Q {\mathbf{B}},\quad G = 2\mathbf A^\top\bf Q \mathbf {B},
\]
\[
L = \mathbf{F}+\mathbf{E}\mathbf{B},\quad M = \mathbf{E}\mathbf{A}, \quad c^\top = \underbrace{[b^\top,\ldots,b^\top]}_{N_p+1\; \mr{times}}\,,
\]
\[
\mathbf{A} = 
\begin{bmatrix}
I \\ A \\ A^2\\ \vdots\\ A^{N_p}
\end{bmatrix},\qquad 
\mathbf{B} = \begin{bmatrix}0&0&\ldots&0\\
B &0&\ldots&0\\
AB & B & \ldots  & 0 \\
 \vdots &\ddots&\ddots\\
 A^{N_p-1}B & \ldots & AB&B \end{bmatrix},
\]
\[
\mathbf{Q} = I_{N_p+1}\otimes Q,\;\; \mathbf{R} = I_{N_{p}}\otimes R,
\]
\[
\mathbf{E} = I_{N_p+1}\otimes E,\;\; \mathbf{F} = \begin{bmatrix}I_{N_{p}}\otimes F \\ 0_{2(N+m)\times mN_p}\end{bmatrix},
\]
\[
F = \begin{bmatrix}0_{2N\times m}\\I_m\\-I_m\end{bmatrix},\;\; E = \begin{bmatrix}I_N\\-I_N\\0_{2m\times N}\end{bmatrix},\;\; b = \begin{bmatrix}
z_{\mr{max}}\\-z_{\mr{min}}\\u_{\mr{max}}\\-u_{\mr{min}}
\end{bmatrix}
\]
with $\otimes$ denoting the Kronecker product (i.e., $I_N\otimes \cdot$ is the N-fold block-diagonalization operator). This form is particularly suited for the active set solver qpOASES (\cite{ferreau2014qpoases}), which also allows for efficient warm-starting.

\subsection{Predictor construction}\label{sec:predConst}
In this section we describe how to construct the predictor of the form~(\ref{eq:pred_lin}) from measured data. We assume that data of the form
 \begin{equation}\label{eq:data}
 \bf X = \begin{bmatrix} x_1,\;\ldots, x_K \end{bmatrix},\; \bf Y = \begin{bmatrix} y_1,\ldots, y_K \end{bmatrix},\; \bf U = \begin{bmatrix} u_1,\ldots, u_K \end{bmatrix}
\end{equation}
where $x_i = [\delta_i^\top,\omega_i^\top]^\top$ and $y_{i} = f(x_i,u_i)$ with $f$ being the dynamics discretized with sampling period $T_{\rm s}$. In other words, if $x_i = [\delta(t)^\top,\omega(t)^\top]$, then $y_i = [\delta(t+T_s)^\top,\omega(t+T_s)^\top]^\top$, i.e, $(x_i,y_i)$ is a pair of successive measurements of the state produced by the continuous time dynamics with the control input $u_i$ held constant during the sampling period $T_{\rm s}$. No relation between $x_i$ and $x_j$, $i\ne j$, is assumed; the data can but is not required to lie on a single trajectory. The data can be collected from the real system or artificially simulated from a model, if available.

The matrices $A$, $B$ of the predictor~(\ref{eq:pred_lin}) are obtained as the solution to the least-squares problem
\begin{equation}\label{eq:ls_AB}
\min_{A,B} \| {\bf Y}_{\mr{lift}} - A {\bf X}_\mr{lift} - B{ \bf U} \|_F,
\end{equation}
where $\|\cdot\|_F$ denotes the Frobenius norm\footnote{The Frobenius norm of a matrix $A$ is given by  $\|A\|_F = \sqrt{\sum_{i,j} A_{i,j}^2}$.} of a matrix and
\begin{equation}\label{eq:liftMat}
\bf X_{\mr{lift}} = \begin{bmatrix} \bs\psi(x_1),\;\ldots, \bs\psi(x_K) \end{bmatrix},\; \bf Y_{\mr{lift}} = \begin{bmatrix} \bs\psi(y_1),\ldots, \bs\psi(y_K) \end{bmatrix}
\end{equation}
with $\bs\psi$  defined in~(\ref{eq:embedMap}). The matrix $C$ is obtained as
\begin{equation}\label{eq:ls_C}
\min_C \| {\bf X} - C {\bf X}_\mr{lift} \|_F.
\end{equation}
The analytic solution to these least-squares problems is
\begin{equation}\label{eq:ls_anal}
[A, B] = \bf Y_\mr{lift} [\bf X_{\mr{lift}}, \bf U]^\dagger,\quad C = \bf X \bf X_{\mr{lift}}^\dagger,
\end{equation}
where $\cdot^\dagger$ denotes the Moore-Penrose pseudoinverse of a matrix. See~\cite[Section 4.1]{korda2016linear} for a more computationally efficient way to obtain the solution for a large $K$.

\subsection{Distributed Koopman MPC}\label{sec:distMPC}
In order to reduce communication and adhere to privacy requirements, it is natural to consider distributed control where the control inputs are determined by several controllers with only partial information availability for each. In this work we consider two natural possibilities:
\begin{enumerate}
\item One controller per grid with information only from within the same grid.
\item Only the grid where the fault occurs is controlled with information only from within the same grid.
\end{enumerate}


\section{Numerical results}
This section summarizes numerical results of our experiments. In order to design the distributed Koopman MPC controllers we construct the predictors of the form (\ref{eq:pred_lin}) from data as descried in Section~\ref{sec:predConst}. The data set consists of $10^4$ trajectories of length $2.5\,\mr{s}$ sampled with period $T_s = 50\,\mr{ms}$ which were collected using the model in the pre-fault configuration. The initial conditions of the trajectories are drawn  uniformly at random with each $\delta \in [-\pi/10,\pi/10]$ and $\omega \in [-0.05,0.05] $. The control input $u_{ij}$ in~(\ref{eqn-SSK:SE_7NE}) is constrained to $[-0.2,0.2]$, i.e., we allow at most $20\,\%$ deviation from the nominal mechanical input power. When constructing the predictors, each control input is distributed uniformly at random withing its bounds.

\begin{figure}[htb]
	\begin{picture}(50,450)
	\put(-70,-20){\includegraphics[width=220mm]{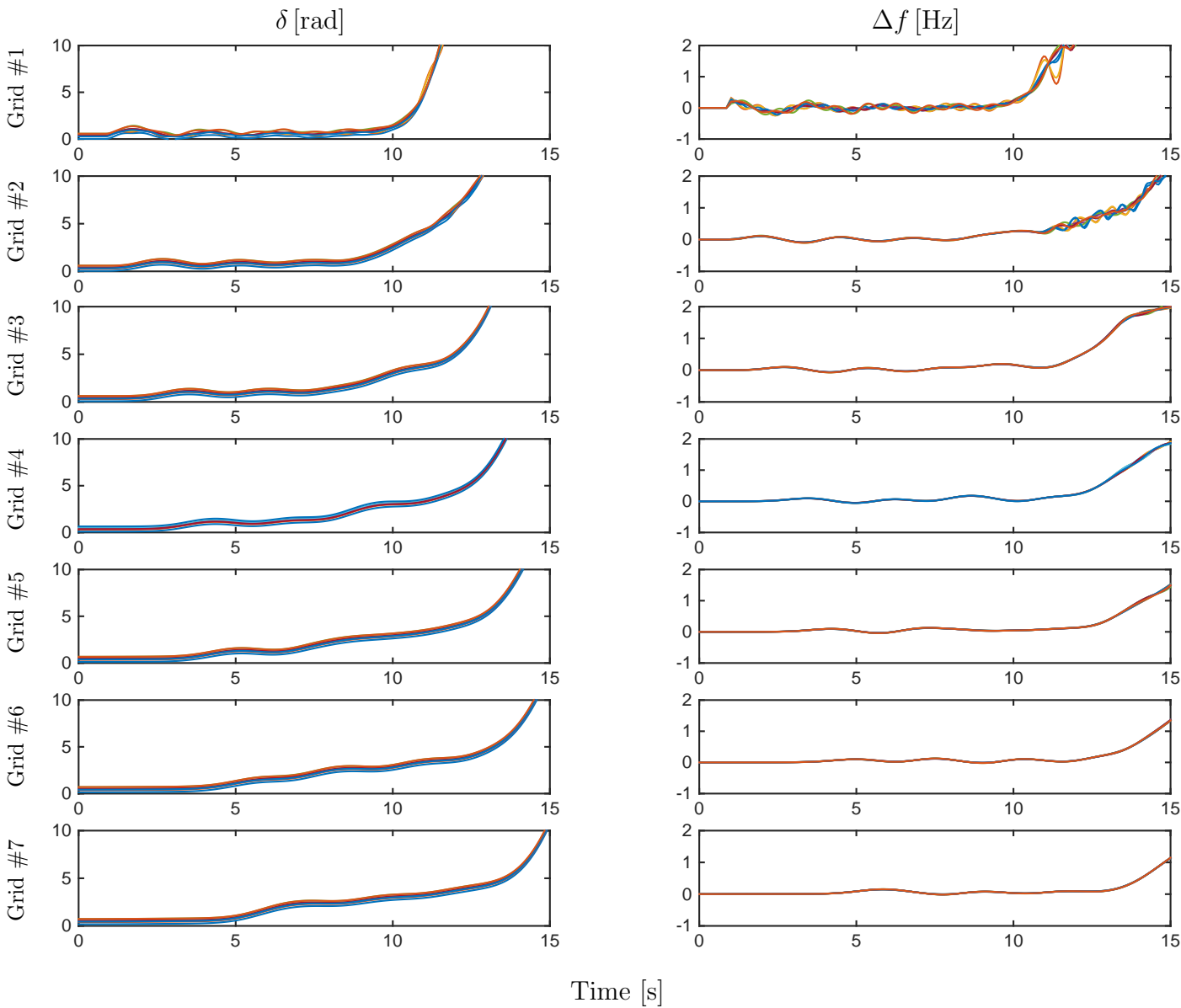}}
	\put(230,0){Time [s]}
	\put(-20,390){\rotatebox{90}{\footnotesize Grid $\#1$}}
	\put(-20,330){\rotatebox{90}{\footnotesize Grid $\#2$}}
	\put(-20,270){\rotatebox{90}{\footnotesize Grid $\#3$}}
	\put(-20,210){\rotatebox{90}{\footnotesize Grid $\#4$}}
	\put(-20,155){\rotatebox{90}{\footnotesize Grid $\#5$}}
	\put(-20,95){\rotatebox{90}{\footnotesize Grid $\#6$}}
	\put(-20,35){\rotatebox{90}{\footnotesize Grid $\#7$}}
	
	\put(95,435){ $\delta\,$[rad]}
	\put(360,435){ $\Delta f\,$[Hz]}
	\end{picture}
	\caption{\footnotesize Cascade instability in the interconnection of seven New England power grids from Figure~\ref{fig:NEcascade}.}
	\label{fig:noCont}
\end{figure}

\begin{figure}[htb]
	\begin{picture}(50,400)
	\put(-70,-20){\includegraphics[width=220mm]{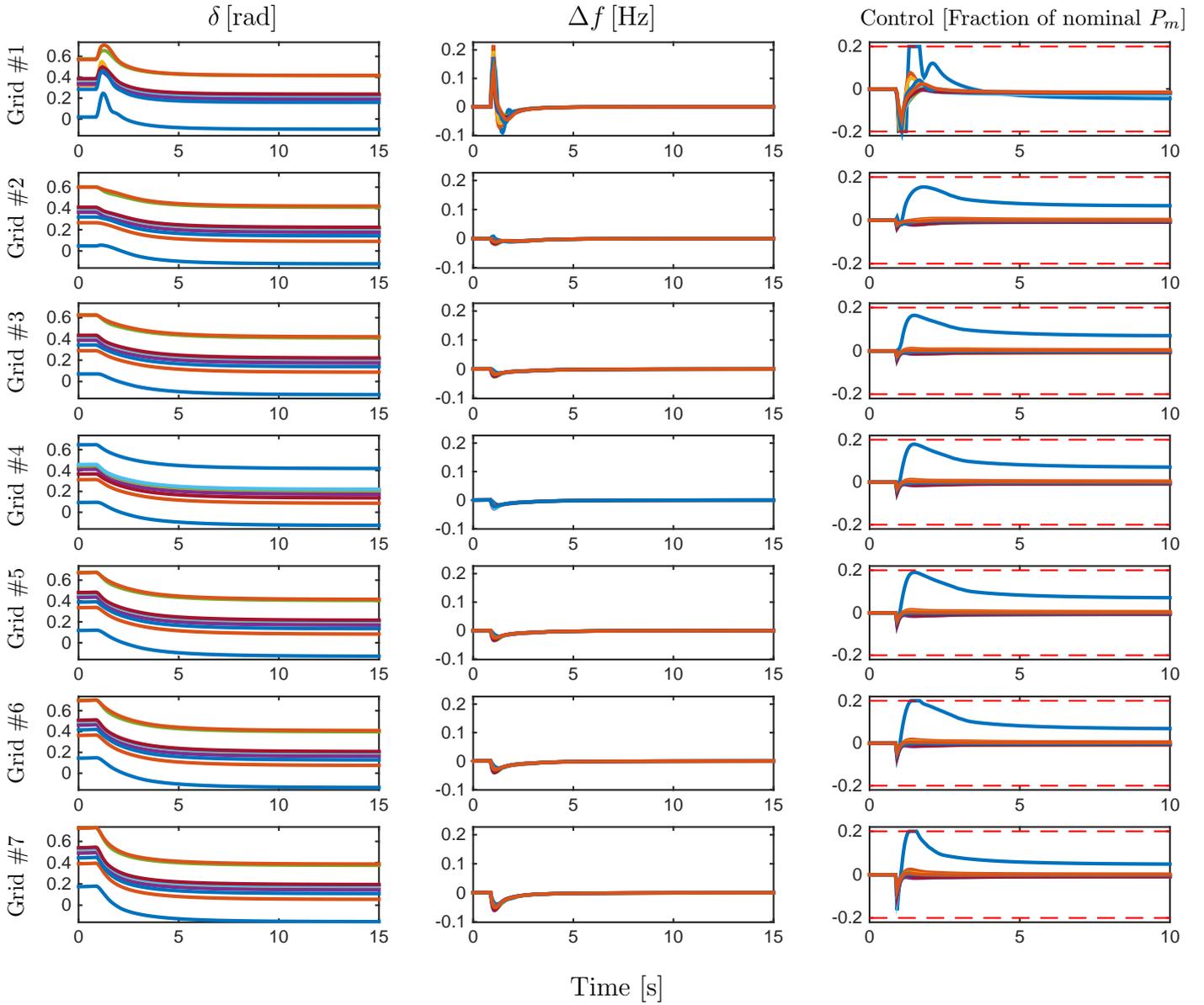}}
	\put(230,0){Time [s]}
	\put(-20,390){\rotatebox{90}{\footnotesize Grid $\#1$}}
	\put(-20,330){\rotatebox{90}{\footnotesize Grid $\#2$}}
	\put(-20,270){\rotatebox{90}{\footnotesize Grid $\#3$}}
	\put(-20,210){\rotatebox{90}{\footnotesize Grid $\#4$}}
	\put(-20,155){\rotatebox{90}{\footnotesize Grid $\#5$}}
	\put(-20,95){\rotatebox{90}{\footnotesize Grid $\#6$}}
	\put(-20,35){\rotatebox{90}{\footnotesize Grid $\#7$}}
	
	\put(65,435){ $\delta\,$[rad]}
	\put(225,435){ $\Delta f\,$[Hz]}
	\put(355,435){ \footnotesize Control [Fraction of nominal $P_m$]}
	\end{picture}
	\caption{\footnotesize Power grid transient stabilization: distributed Koopman MPC with one controller per grid.}
	\label{fig:ContAllGrids}
\end{figure}

\paragraph{One controller per grid} First, we consider the distributed Koopman MPC with one controller per grid. The cost matrices $Q$ and $R$ are chosen to 
\[
Q = \mr{bdiag}(0_{ n_{\mr{gen}} \times n_{\mr{gen}}  }, 0_{ n_{\mr{gen}}\times n_{\mr{gen}}}, I_{ n_{\mr{gen}}  }), \quad R = 0.01  I_{ n_{\mr{gen}}  },
\]
where $n_{\mr{gen}} = 9$ is the number of generators in each grid, i.e., we penalize $\sum_{i=0}^{N_p-1} \|\omega_i\|_2^2 + 0.1 \|u\|_2^2$.  The prediction horizon is chosen to be $1\,\mr{s}$, corresponding to $N_p = 1 / T_s = 20$. We do not impose any constraints on the state variables whereas the control input, i.e., the mechanical input power of the generators, is constrained to be between $\pm 20\,\%$ of the nominal value for each generator. Simulation results with no control are depicted in Figure~\ref{fig:noCont}; in accordance with~\cite{susuki2012coherent}, we observe unstable behavior caused by the fault at time $t_f = 0.87\,\mr{s}$. Results with control are depicted in Figure~\ref{fig:ContAllGrids}. We observe a stable behavior, very fast attenuation of the disturbance and a bounded maximum frequency deviation ($\approx 0.2\,\mr{Hz}$). The quadratic program~(\ref{opt:MPC}) was solved using qpOASES~\cite{ferreau2014qpoases} running on Matlab. The average computation time required for evaluation of the control input was approximately $10\,\mr{ms}$ on a laptop with macOS, 2GHz intel i7, which would allow for a real-time implementation.

\begin{figure}[htb]
	\begin{picture}(50,400)
	\put(-70,-20){\includegraphics[width=220mm]{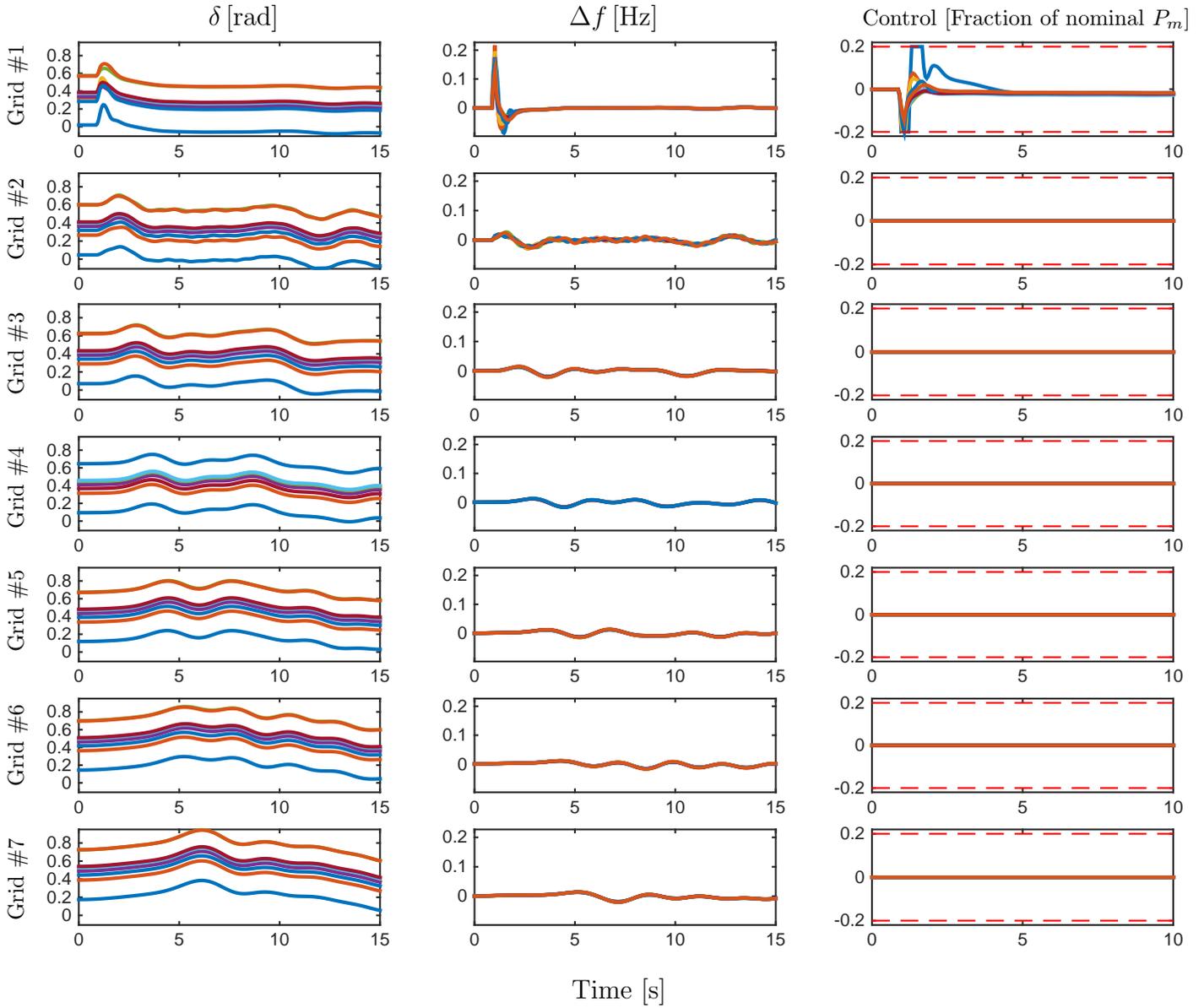}}
	\put(230,0){Time [s]}
	\put(-20,390){\rotatebox{90}{\footnotesize Grid $\#1$}}
	\put(-20,330){\rotatebox{90}{\footnotesize Grid $\#2$}}
	\put(-20,270){\rotatebox{90}{\footnotesize Grid $\#3$}}
	\put(-20,210){\rotatebox{90}{\footnotesize Grid $\#4$}}
	\put(-20,155){\rotatebox{90}{\footnotesize Grid $\#5$}}
	\put(-20,95){\rotatebox{90}{\footnotesize Grid $\#6$}}
	\put(-20,35){\rotatebox{90}{\footnotesize Grid $\#7$}}
	
	\put(65,435){ $\delta\,$[rad]}
	\put(225,435){ $\Delta f\,$[Hz]}
	\put(355,435){ \footnotesize Control [Fraction of nominal $P_m$]}
	\end{picture}
	\caption{\footnotesize Power grid transient stabilization: only first grid controlled by Koopman MPC.}
	\label{fig:ContFirstGrid}
\end{figure}

\paragraph{Only affected grid controlled} Next we consider the scenario when only the  grid where the fault occurred (i.e., the first grid in our case) is controlled and there is no information transmitted among the grids; we use exactly the same controller for the first grid as in the previous case whereas the remaining grids remain uncontrolled. The computation time to evaluate the control input in the first grid remains approximately the same ($\approx 10\,\mr{ms}$), whereas no control action is taken in the remaining grids.  Figure~\ref{fig:ContFirstGrid} shows the results. We observe stable behavior but  the oscillations in the frequency deviations in the uncontrolled grids $\#2-\#7$ persist much longer even though the maximum frequency deviation is of similar magnitude (since the maximum deviation is attained in the first grid which is controlled in both cases).


\clearpage
\bibliographystyle{abbrv}
\bibliography{References}

\end{document}